\documentclass[12pt]{article}
\usepackage{theorem}
\newtheorem{lemma}{Lemma}

\newtheorem{theorem}[lemma]{Theorem}
\newtheorem{corollary}[lemma]{Corollary}
{\theorembodyfont{\upshape}}
{\theorembodyfont{\upshape}}
{\theorembodyfont{\upshape}\newtheorem{conjecture}[lemma]{Conjecture}}
{\theorembodyfont{\upshape}\newtheorem{example}[lemma]{Example}}
{\theorembodyfont{\upshape}}
{\theorembodyfont{\upshape}\newtheorem{note}[lemma]{Note}}
\newcommand{\Z}{{\bf Z}}
\newcommand{\R}{{\bf R}}
\newcommand{\C}{{\bf C}}
\newcommand{\rme}{{\rm e}}
\newcommand{\rmd}{\,{\rm d}}
\newcommand{\cC}{{\cal C}}
\newcommand{\cD}{{\cal D}}

\newcommand{\cH}{{\cal H}}

\newcommand{\sig}{\sigma}
\newcommand{\alp}{\alpha}
\newcommand{\bet}{\beta}
\newcommand{\gam}{\gamma}
\newcommand{\lam}{\lambda}
\newcommand{\del}{\delta}
\newcommand{\eps}{\varepsilon}
\newcommand{\Gam}{\Gamma}
\newcommand{\kap}{\kappa}

\newcommand{\Dom}{{\rm Dom}}

\newcommand{\Spec}{{\rm Spec}}

\newcommand{\norm}{\Vert}
\renewcommand{\Re}{{\rm Re}\;}

\newcommand{\Proof}{\underbar{Proof}{\hskip 0.1in}}

\newcommand{\Schrodinger}{Schr\"odinger }
\newcommand{\pr}{\prime}
\newcommand{\emp}[1]{{\it #1}}
\newcommand{\plus}{_+}
\newcommand{\minus}{_-}
\title{EIGENVALUES OF \\ AN ELLIPTIC SYSTEM}
\author{E.B. Davies}
\date{August 2001}
\begin{document}
\maketitle
%%%%%%%%%%%%%%%%%%%%%
\begin{abstract}
We describe the spectrum of a non-self-adjoint elliptic system on a
finite interval. Under certain conditions we find that the
eigenvalues form a discrete set and converge asymptotically at
infinity to one of several straight lines. The eigenfunctions need
not generate a basis of the relevant Hilbert space, and the larger
eigenvalues are extremely sensitive to small perturbations of the
operator. We show that the leading term in the spectral asymptotics
is closely related to a certain convex polygon, and that the
spectrum does not determine the operator up to similarity. Two
elliptic systems which only differ in their boundary conditions may
have entirely different spectral asymptotics. While our study makes
no claim to generality, the results obtained will have to be
incorporated into any future general theory.
\par
keywords: eigenvalues, non-self-adjoint operators, ordinary
differential operators, spectral theory, basis, pseudospectrum.
\end{abstract}
%%%%%%%%%%%%%%%%%%%%%%%%
\section{Introduction}
\par
Consider a first order elliptic differential system acting in
$L^2:=L^2((\alp,\bet),\C^n)$ according to the formula
\[
(Lf)(x):=A(x)f^\pr(x)
\]
where $A$ is an $n\times n$ complex matrix depending on $x$. (We
show in Section 2 that our methods also apply to certain second and
higher order elliptic systems.) If $A(x)$ is a piecewise continuous
function of $x$ then the solution to the equation $Lf=zf$ can be
written in the form
\[
f(x)=U(z,x)f(\alp)
\]
where $U(z,x)$ is an invertible $n\times n$ matrix which depends
continuously on $x$ and analytically on $z$. The operator has
domain contained in $W^{1,2}:=W^{1,2}((\alp,\bet),\C^n)$. Since
$W^{1,2}\subseteq C[\alp,\bet]$, we can impose general boundary
conditions of the type
\begin{equation}
Sf(\alp)+Tf(\bet)=0. \label{bc}
\end{equation}
where $S,T$ are any $n\times n$ matrices.

Problems of this type arise in several contexts. One relates to
non-equilibrium thermodynamics, \cite{str} and another concerns the
study of generalized determinants for elliptic operators on
manifolds with boundary, \cite{lesch,scott,scott2}. Of course
problems involving higher order differential operators can also be
reformulated in terms of first order systems; we do not, however,
assume Hamiltonian structure, which would force $n$ to be even and
also imply strong constraints on the solutions of the equations.
The present paper provides a class of nearly exactly soluble
examples which illustrate some of the phenomena such studies have
to face.

The spectral behaviour of $L$ has obvious relevance to the
time-dependent system
\begin{equation}
\frac{\partial f}{\partial t}=A(x)\frac{\partial f}{\partial
x}\label{hyperb}.
\end{equation}
If each matrix $A(x)$ has only real eigenvalues $a_r(x)$ then
(\ref{hyperb}) is called a hyperbolic system and $a_r(x)$ are the
characteristic speeds at $x$. This case has special features by
comparison with the case in which the eigenvalues are in general
position in the complex plane, and is only partially analyzed in
this paper. See Section 7 on semigroup properties.

\begin{example}\label{subs}
Consider the `Dirichlet' boundary conditions $f(\alp)\in U$,
$f(\bet)\in V$, where $U,\, V$ are linear subspaces in $\C^n$ with
\[
\dim (U)+\dim (V)=n.
\]
This falls within the above scheme if we choose $S,T$ as follows:
$S$ should have kernel $U$ and range $W$ while $T$ should have
kernel $V$ and range $X$ where $W\cap X=\{0\}$ and $W+X=\C^n$. In
particular if $U$ is the linear span of the single vector $u\not=
0$ and $V$ is the annihilator of the vector $v\not= 0$ then the
boundary value problem reduces to
\[
F(z):=\langle U(z,\bet )u,v\rangle =0.
\]
\end{example}

Returning to the general context, the embedding of $W^{1,2}$ into
$L^2$ is compact, so the operator $L$ has compact resolvent if its
spectrum is not equal to $\C$, and its spectrum must consist of
discrete eigenvalues of finite multiplicity. If $A(x)=a I$ for all
$x\in [\alp,\bet]$ then the only possible solution of $Lf=z f$ is
$f(x)=c\rme^{z x/a}$ for some $c\in
\C^n$; in the context of Example \ref{subs}
we deduce that $\Spec(L)=\C$ if $U\cap V\not=\{ 0\}$, but in
all other cases the spectrum is empty. The situation changes
entirely if we assume that $A(x)$ is constant but its eigenvalues
$a_r$ are all different and all non-zero.

It follows immediately from the above discussion that $z$ is an
eigenvalue of $L$ if and only if
\[
F(z):=\det\left(S+TU(z,\bet )\right)=0
\]
where $F$ is an entire function of $z$. There are well-developed
numerical procedures for evaluating $F(z)$ for any $z$ and for
finding the points at which $F$ vanishes. However, the asymptotic
behaviour of the eigenvalues of $L$ and the study of such
quantities as the regularized determinant are difficult subjects.
There is one case in which they can be analyzed fairly directly,
namely when $A(x)$ is piecewise constant, an assumption which we
make for the remainder of the paper.

In the following sections we present a detailed analysis of the
function $F(z)$. We prove in the `generic case' that there are
several series of eigenvalues diverging to infinity along straight
lines which are determined by a certain convex polygon. The
asymptotic form of the spectral counting function depends upon the
length of the boundary of this polygon. There are exceptional cases
in which almost periodic structure arises, and which we do not
analyze fully. The implications of our results for the spectral
analysis are explained with a series of simple examples.

A much more complete analysis of the second order case with $n=2$
has been performed by Boulton, \cite{bou,bou2}, and his results
have informed our approach to the case $n\geq 2$. Our results
describe how the more general analysis goes in the generic case,
and what further phenomena have to be considered for a complete
analysis.

\section{The function $F(z)$}

We assume henceforth that
\[
\alp=\alp_0<\alp_1<...<\alp_m=\bet
\]
and that $A(x)=A_s$ if $\alp_{s-1}<x\leq \alp_s$. We also assume
that each $A_s$ is invertible and diagonalizable. Then
\[
U(z,\bet )=U_mU_{m-1}...U_2U_1
\]
where
\[
U_s=\rme^{zA_s^{-1}(\alp_s-\alp_{s-1})}=V_s\rme^{zD_s}V_s^{-1}
\]
and each matrix $D_s$ is diagonal. We then have
\begin{eqnarray}
F(z)&=&\det\left(S+TV_m\rme^{zD_m}V_m^{-1}.V_{m-1}
\rme^{zD_{m-1}}V_{m-1}^{-1}...V_1\rme^{zD_1}V_1^{-1}
\right)\nonumber\\
&=&\sum_{r=1}^R\del_r\rme^{z\overline{\gam_r}}.\label{piecewise}
\end{eqnarray}

It might be thought that the case in which $A(x)$ is continuous can
be obtained from the above by an approximation procedure. In
Example \ref{pathol} we show that taking such a limit would not be
a straightforward matter. We discuss the spectral asymptotics of
such piecewise constant operators in Section 4, and make a
conjecture about the case in which $A(x)$ depends continuously on
$x$.

Let us look at the case of Dirichlet boundary conditions in more
detail. Assuming that $A$ is independent of $x$, invertible and
diagonal with eigenvalues $a_r$, the solution of the eigenvalue
equation is
\[
f_r(x)=c_r\rme^{z(x-\alp )/a_r}
\]
where $c\in \C^n$ must be chosen so that $f$ satisfies the boundary
conditions. If $u_1,...,u_{n-p}$ is a basis of the annihilator of
$U$ in the dual space and $v_1,...,v_{p}$ is a basis of the
annihilator of $V$ then the conditions are that $(c, u_i) =0$ for
all $1\leq i \leq n-p$ and $(c, w_j)=0$ for all $1\leq j\leq p$,
where $w_j\in\C^n$ is defined by
\[
w_{j,r}=v_{j,r}\rme^{z(\bet-\alp )/a_r}.
\]
The existence of a non-zero eigenvector then reduces to
\[
\det\{ v_1,...,v_{n-p},w_1,...,w_{p}\}=0.
\]
This may be rewritten in the form $F(z)=0$ where $F$ is defined by
(\ref{Fz}). In this case each of the coefficients $\gam_r$ is a sum
of $p$
\emp{distinct} terms of the form $(\bet-\alp )/\overline{a_i}$.

In the simplest case $\dim (U) =1$ and there are $n$ distinct
$\gam_r$. In general the number of $\gam_r$ is given by a
combinatorial expression. However, some of the $\gam_r$ may
coincide, and some of the $\del_r$ may vanish.

We next consider a similar problem for second order elliptic
systems. Let $H$ act in $L^2:=L^2((\alp,\bet),\C^n)$ according to
the formula
\[
(Hf)_r(x):=a_r^2\frac{\rmd^2f_r}{\rmd x^2}
\]
for $1\leq r\leq n$, where $a_r\not= 0$ for all $r$.  The operator
has domain contained in $W^{2,2}:=W^{2,2}((\alp,\bet),\C^n)$. Since
$W^{2,2}\subseteq C^{(1)}[\alp,\bet]$, the boundary conditions may
involve $f(\alp ),\, f^\pr (\alp ),\, f(\bet ),\, f^\pr (\bet )$.
We do not consider the most general choice of boundary conditions,
which leads to a somewhat complicated function $F(z)$, but follow
Boulton in assuming a combination of Dirichlet or Neumann boundary
conditions in the following sense, \cite{bou,bou2}. We suppose that
$f(\alp )\in U_1$, $f^\pr(\alp )\in U_2$, $f(\bet )\in V_1$ and
$f^\pr(\bet )\in V_2$ where $U_1,\, U_2,\, V_1,\, V_2$ are linear
subspaces  of $\C^n$. We impose the minimal further condition
\[
\dim(U_1)+\dim(U_2)+\dim(V_1)+\dim(V_2)=2n.
\]

The eigenvalue problem may be written in the form $Hf=z^2f$, and
its solution is of the form
\[
f_r(x)=c_r\rme^{(x-\alp )z/a_r}+d_r\rme^{-(x-\alp )z/a_r}
\]
where $(c,d)\in \C^{2n}$, provided $z\not= 0$. The case $z=0$ is
dealt with separately. The boundary conditions lead to $2n$ linear
constraints on the vector $(c,d)$ which have a non-zero solution if
$F(z)=0$, where $F(z)$ is a function of the form (\ref{Fz}).

As a particular case we mention the possibility $H:=L^2$ where $L$
is the first order operator already discussed. This corresponds to
the choice of boundary conditions $f(\alp )\in U$, $f^\pr(\alp )\in
A^{-1}U$, $f(\bet )\in V$, $f^\pr (\bet )\in A^{-1}V$, where
$A_{i,j}=\del_{i,j}a_i$. The eigenvalues of $H$ are the squares of
the eigenvalues of $L$, and it follows from our analysis below that
they are asymptotic to certain parabolae at infinity.

\section{Zeros of $F(z)$}

Our task is to study the asymptotic distribution of the zeros of
entire functions of the form
\begin{equation}
F(z)=\sum_{r=1}^R\del_r\rme^{z\overline{\gam_r} }=\sum_{r=1}^R
F_r(z) \label{Fz}
\end{equation}
as $|z|\to\infty$. We assume that all $\gam_r$ are different, that
all $\del_r$ are non-zero and that $R\geq 2$. We will see that the
asymptotic structure of the set of zeros depends heavily on the
convex hull $K$ of the set of $\gam_r$. We will only consider the
\emp{generic case}, defined as follows. After relabelling, the
convex hull has vertices $\gam_1,...,\gam_Q$ written in
anticlockwise order, while the remainder of the $\gam_r$ lie in the
interior of $K$; we assume that no $\gam_r$ lies on a edge of $K$
unless it is already a vertex of $K$. If $1\leq r\leq Q$ we define
\begin{eqnarray*}
r\plus &:=&\left\{\begin{array}{ll} r+1&\mbox{if $1\leq r\leq
Q-1$}\\ 1&\mbox{if $r=Q$,}
\end{array}\right.\\
r\minus &:=&\left\{\begin{array}{ll} r-1&\mbox{if $2\leq r\leq
Q$}\\ Q&\mbox{if $r=1$}.
\end{array}\right.
\end{eqnarray*}
We will prove that, in addition to a finite number of zeros near to
the origin, $F(z)$ has $Q$ series of zeros, each of which converges
asymptotically at infinity towards one of the sets on which
\[
|F_r(z)|=|F_{r\minus}(z)|.
\]
where $1\leq r\leq Q$. Rewriting this in the form
\[
| \del_r|\rme^{z\cdot\gam_r } =|
\del_{r\minus}|\rme^{z\cdot\gam_{r\minus}}
\]
we see that it is a straight line
\begin{equation}
z\cdot (\gam_r-\gam_{r\minus}) =k_r  \label{line}
\end{equation}
perpendicular to the edge $(\gam_{r\minus},\gam_r)$ of $K$.

If one or more of the $\gam_r$ with $Q<r\leq R$ lie within an edge
of $K$ then the spectral behaviour of $L$ at infinity is much more
complicated. Boulton has shown that the behaviour depends upon
whether those $\gam_r$ which lie within an edge of $K$ divide that
edge into parts whose lengths have rational or irrational ratios.
In the former case the eigenvalues still converge at infinity to
one of several straight lines, but there are several lines
perpendicular to the relevant edge instead of just one. In the
irrational case the eigenvalues have almost periodic structure at
infinity within one of several strips, with one strip perpendicular
to each edge. Boulton's analysis is only presented in a typical
second order case, but the ideas clearly extend as stated,
\cite{bou,bou2}.

Returning to the generic case, we study the asymptotic location of
the zeros of $F(z)$ by dividing the exterior of $K$ into $4Q$
subregions. Given $\eps
>0$ there are $3$ semi-infinite strips $S_r^\eps,\, S_r^+,\,
S_r^-$ associated with each edge $(\gam_{r\minus},\gam_r)$ of $K$,
and also a wedge $W_r$ associated with each vertex $\gam_r$. Let
$e_r$ be the outward pointing unit normal to the edge
$(\gam_{r\minus},\gam_r)$. Given $1\leq r\leq Q$ we define
\begin{eqnarray*}
S_r^\eps&:=& \{ z\in\C:z\cdot e_r\geq \gam_r\cdot e_r \mbox{ and }
k_r-\eps \leq z\cdot (\gam_r-\gam_{r\minus})\leq k_r+\eps\}\\
S_r^+&:=& \{ z\in\C:z\cdot e_r\geq \gam_r\cdot e_r \mbox{ and }
k_r+\eps
\leq z\cdot (\gam_r-\gam_{r\minus})\leq k_r+1\}\\
S_r^-&:=&
\{ z\in\C:z\cdot e_r\geq \gam_r\cdot e_r \mbox{ and } k_r-1 \leq
z\cdot (\gam_r-\gam_{r\minus})\leq k_r-\eps\}\\ W_r&:=&
\{z\notin K:k_r+1\leq z\cdot (\gam_r-\gam_{r\minus})
\mbox{ and } k_{r\plus}-1\geq z\cdot (\gam_{r\plus}-\gam_r)\}.
\end{eqnarray*}

\begin{theorem}
Given $0<\eps<1/9$ there exists $T_\eps$ such that any eigenvalue
$z$ of $L$ satisfying $|z|>T_\eps$ must lie in $S_r^\eps$ for some
$r$ such that $1\leq r\leq Q$.
\end{theorem}

\Proof If $z\in S_r^+$ then
\[
|F_{r\minus}(z)/F_{r}(z)|\leq\rme^{-\eps}<1-\eps/2.
\]
If $s\not= r,r\minus$ then $(\gam_r-\gam_s)\cdot e_r>0$ and for
large enough $|z|$ we have
\[
|F_{s}(z)/F_{r}(z)|\leq\eps/2R.
\]
Combining these bounds we see that for such $z$
\[
|F(z)|\geq |F_r(z)|\eps(1/2-(R-2)/2R)>0.
\]
A similar argument applies to $S_r^-$.

If $s\not= r,r_\pm$ then either
\[
\gam_r-\gam_s=p(\gam_r-\gam_{r\minus})+q e_r
\]
where $p>0$ and $q >0$, or
\[
\gam_r-\gam_s=p(\gam_r-\gam_{r\plus})+q e_{r\plus}
\]
where $p>0$ and $q >0$. It follows that
\[
\lim_{|z|\to\infty,\,z\in W_r} \left\{
(\gam_r-\gam_s)\cdot z\right\}=+\infty.
\] Now let $z\in W_r$. We have
\[
|F_{r_\pm}(z)/F_r(z)|\leq \rme^{-1}
\]
and if $|z|$ is large enough we also have
\[
|F_{s}(z)/F_{r}(z)|\leq\eps/2R.
\]
provided $s\not= r,r_\pm$. Combining these inequalities we obtain
\[
|F(z)|\geq |F_r(z)|(1-2\rme^{-1}-\eps(R-3)/2R)>0.
\]

\begin{note}
If one of the $\gam_s$ for $Q<s\leq R$ is very close to
the boundary of $K$ then $T_\eps$ will be correspondingly large in
the above proof. This is inevitable because the possible asymptotic
directions of the spectrum change as $\gam_s$ moves through an edge
of $K$, whereupon $K$ changes so that $\gam_s$ becomes another
vertex of $K$.
\end{note}

\begin{note}
The semi-infinite lines (\ref{line}) divide the complex plane into $Q$
sectors $S_r$. If $|z|\to \infty$ within $S_r$ then $F(z)\sim
\del_r\rme^{z\overline{\gam_r}}$, so
\[
\log(F(z))\sim \log(\del_r)+z\overline{\gam_r}.
\]
This formula may be used to compute the regularized determinant in
the sense of \cite{lesch,scott,scott2}.
\end{note}

\section{Spectral Asymptotics}

We have proved that the zeros of $F(z)$ converge asymptotically
towards one of $Q$ straight lines as $|z|\to \infty$. To obtain
more detailed information we introduce the functions
\[
G_r(z):=F_r(z)+F_{r\minus}(z).
\]

\begin{theorem}\label{gr}
The zeros of  $G_r(z)$ lie on the line (\ref{line}) with constant
distance $2\pi/|\gam_r-\gam_{r\minus}|$ between any two consecutive
zeros.
\end{theorem}

\Proof If we put $\gam_r-\gam_{r\minus}=\rho\rme^{i\theta}$ where
$\rho >0$ and $\theta\in\R$ and $c_r=\log(\del_{r\minus}/\del_r)$,
then the zeros of $G_r(z)$ are given by
\[
z=(c_r+(2n+1)\pi i)\rme^{i\theta}/\rho
\]
where $n\in\Z$. The statements of the lemma all follow from this.

\begin{theorem}\label{f}
The zeros of $F(z)$ which lie in $S_r^\eps$ converge as
$|z|\to\infty$ to the zeros of $G_r(z)$, in the sense that the
modulus of the differences of corresponding zeros converges to
zero.
\end{theorem}

\Proof We have already proved that the zeros of $F(z)$ converge to
one of the lines (\ref{line}). The statement is a straightforward
application of Rouche's theorem, since we have already noted that
the remaining terms in the series are asymptotically negligible as
$|z|\to\infty$ within $S_r^\eps$.

\begin{corollary}\label{asy}
If $N(E)$ is the number of zeros of $F(z)$ such that $|z|\leq E$
then
\[
N(E)=b(K)E/2\pi +O(1)
\]
as $E\to\infty$, where $b(K)$ is the length of the boundary of $K$.
\end{corollary}

\Proof The number of zeros of $F(z)$  satisfying $|z|\leq E$ associated with
each line is $|\gam_r-\gam_{r\minus}|E/2\pi +O(1)$ by the above
theorem, since the zeros of $F(z)$ only converge to the zeros of
$G_r(z)$ in one direction. The corollary follows by summing over
$r$.
\begin{example}
If $L_i$ are two elliptic systems acting in
$L^2((\alp,\bet),\C^{n_i})$ respectively for  $i=1,2$ then we may
consider their direct sum $L=L_1\oplus L_2$. By applying Corollary
\ref{asy} to each of the components we obtain
\[
N(E)=(b(K_1)+b(K_2))E/2\pi +O(1)
\]
in an obvious notation. On the other hand we have
$F(z)=F_1(z)F_2(z)$ so the exponents $\gam_r$ in $F(z)$ are the
sums of the exponents in $F_1(z)$ and in $F_2(z)$. This implies
that $K=K_1+K_2$. These different approaches to the asymptotics
are reconciled by the classical but non-trivial fact that
\[
b(K_1+K_2)=b(K_1)+b(K_2)
\]
for any two plane convex sets. This example may be used to
construct counterexamples to various conjectures.
\end{example}

\begin{example}\label{pathol}
Let us consider the entire function
\[
F_n(z)=n^{-1}\sum_{r=1}^n\exp\left(\rme^{2\pi ir/n}z\right).
\]
According to Corollary \ref{asy}
\[
N_n(E)=b_nE/2\pi +O(1)
\]
as $E\to\infty$, where $b_n\to 2\pi$ as $n\to\infty$. On the other
hand
\[
\lim_{n\to\infty}F_n(z)=F(z):=\frac{1}{2\pi}\int_0^{2\pi}\exp\left(\rme^{i\theta}z\right)\rmd
\theta.
\]
Since $F(z)$ is rotationally invariant and entire with $F(0)=1$ it must be
identically equal to $1$, and so it has no zeros at all. This
establishes that even for this class of entire functions, the asymptotic
distribution of the zeros does not vary continuously with the
function.
\end{example}

When combined with Corollary \ref{asy}, the following theorem
proves that the eigenvalues of $L$ move off to infinity as the
eigenvalues of $A$ coalesce. The singular behaviour of the spectrum
in the limit $\del\to 0$ is in accordance with the behaviour for
$\del =0$ described in the Introduction.

\begin{theorem} Let $A=\alp I+\del B$, where $\alp\not= 0$, the
eigenvalues of $B$ are all distinct, and $\del\in\C$ is
sufficiently small. Let $L$ be defined in the usual manner with
$\dim(U)=p$. Then there exist distinct constants $\sig_r$ such that
\begin{equation}
\gam_r =p(\bet-\alp )/\alp+\sig_r\del+O(\del^2)\label{ga}
\end{equation}
and there exists $\kap >0$ such that
\begin{equation}
b(K)=\kap\del+O(\del^2).\label{epsb}
\end{equation}
\end{theorem}

\Proof The eigenvalues of $A$ are $a_r=\alp+\del b_r$ for
$1\leq r\leq n$ where $b_r$ are the distinct eigenvalues of $B$.
The formula for the $\gam_r$ in Section 2 yields (\ref{ga})
immediately. This implies (\ref{epsb}) with
\[
\kap=\sum_{r=1}^Q|\sig_r-\sig_{r\minus}|/2\pi.
\]

In the remainder of this section we show how to apply Corollary
\ref{asy} to certain elliptic systems with variable coefficients.
We assume that
\[
Lf(x)=A(x)f^\pr(x)
\]
in $L^2((\alp,\bet),\C^n)$, where $A(x)$ is invertible and
diagonalizable for each $x\in (\alp,\bet)$. In our next theorem we
assume boundary conditions of the form $f(\alp)=u$ and $\langle
f(\bet),v\rangle
=0$, where $u,v\in\C^n$ are both non-zero.

The meaning of the word `generically' below will be explained
during the proof.

\begin{theorem}\label{asymptheorem}
Let $b(x)$ be the length of the boundary of $K(x)$, where $K(x)$ is
the convex hull of the eigenvalues of $A(x)^{-1}$. If also
$A(\cdot)$ is piecewise constant then generically one has
\begin{equation}
N(E)=\frac{E}{2\pi}\int_\alp^\bet b(x)\rmd x +O(1)\label{intasymp}
\end{equation}
as $E\to\infty$.
\end{theorem}

\Proof Following the notation of Section 2, equation (\ref{piecewise})
becomes
\begin{eqnarray}
F(z)&=&\langle V_m\rme^{zD_m}V_m^{-1}.V_{m-1}
\rme^{zD_{m-1}}V_{m-1}^{-1}...V_1\rme^{zD_1}V_1^{-1}
u,v \rangle\nonumber\\
&=&\sum_{r=1}^R\del_r\rme^{z\overline{\gam_r}}.\label{newpiece}
\end{eqnarray}
Denoting the eigenvalues of $A_s$ by $\{a_{s,t}\}_{t=1}^n$, $D_s$
is a diagonal matrix with entries $\{
a_{s,t}^{-1}(\alp_s-\alp_{s-1})\}_{t=1}^n$. Each
$\overline{\gam_r}$ in (\ref{newpiece}) is of the form
\begin{equation}
\overline{\gam_r}=\sum_{s=1}^m
a_{s,t(s)}^{-1}(\alp_s-\alp_{s-1})\label{summup}
\end{equation}
where $t(\cdot )$ is a function from $\{1,...,m\}$ to
$\{1,...,n\}$. The parameter $r$ is a relabelling of the set of all
such functions, so $R\leq n^m$, with equality unless two sums of
the form (\ref{summup}) happen to be equal. If $K$ is the convex
hull of the $\gam_r$ then our generic assumption is that
$\del_r\not= 0$ for all vertices $\gam_r$ of $K$, and that no
$\gam_r$ lie within any edge of $K$.

By Corollary \ref{asy} we have to prove that
\begin{equation}
b(K)=\int_\alp^\bet b(x)\rmd x.\label{asympint}
\end{equation}
The right hand side equals
\[
\sum_{s=1}^m b(K_s)(\alp_s-\alp_{s-1})
\]
where $A(x)=A_s$ and $K(x)=K_s$ if $\alp_{s-1}<x\leq \alp_s$ as in
Section 2. The identity (\ref{asympint}) follows by combining the
facts that
\begin{eqnarray}
b(tK_1)&=&tb(K_1)\label{convscale}\\ b(K_1+K_2)&=&
b(K_1)+b(K_2)\label{convsum}
\end{eqnarray}
for any convex sets $K_1,K_2$ and any $t>0$, with the identity
\[
K=\sum_{s=1}^m (\alp_s-\alp_{s-1})K_s.
\]

\begin{conjecture}
We conjecture that Theorem \ref{asymptheorem} remains generically
valid for piecewise continuous coefficients $A(x)$, $x\in
[\alp,\bet]$.
\end{conjecture}

\begin{note}
The leading coefficient in (\ref{intasymp}) depends only on the
symbol of $L$. However, the truth of the theorem depends upon the
genericity assumption for the following reason. Let $L,\tilde{L}$
be two operators with the same symbol but different boundary
conditions. Since the coefficients $\del_r$ in the expansion
(\ref{newpiece}) may vanish for different values of $r$ for the two
operators, it may happen that $K\not= \tilde{K}$ in an obvious
notation. The leading coefficient in the spectral counting function
(\ref{intasymp}) will then usually be different.
\end{note}

This phenomenon is illustrated by the following theorem. We make
the same assumptions as in Theorem \ref{asymptheorem}, except that
we now allow general boundary conditions of the form (\ref{bc}).
The following theorem may be regarded as a version of Weyl's
formula for the operator associated with the matrix-valued symbol
\[
p_{r,s}(x,\xi)=A_{r,s}(x)\xi
\]
where $x\in [\alp,\bet]$ and $\xi\in\R$. Once again the meaning of
the word `generically' will be explained during the proof. We warn
the reader that the asymptotic form obtained here is different from
that of Theorem \ref{asymptheorem}; this is possible because the
boundary conditions of that theorem are all non-generic in the
sense used in the present theorem.
%The self-adjoint (i.e. hyperbolic) case was proved by Ivrii,
%without any of the restrictions imposed in our case; compare
%\cite[Theorem 3]{ds}.

\begin{theorem}\label{genericasymp}
Under the above assumptions one generically has
\[
N(E)=\frac{E}{2\pi}\int_\alp^\bet b(x)\rmd x +O(1)
\]
where
\[
b(x)=2\sum_{r=1}^n\left|a_r(x)\right|^{-1}
\]
and $\{a_r(x)\}_{r=1}^n$ are the eigenvalues of $A(x)$ for each
$x\in(\alp,\bet)$.
\end{theorem}

\Proof
A more detailed understanding of (\ref{piecewise}) may be obtained
by considering
\begin{eqnarray}
\tilde F(z_1,...,z_m)&=&\det\left(S+TV_m\rme^{z_mD_m}V_m^{-1}.V_{m-1}
\rme^{z_{m-1}D_{m-1}}V_{m-1}^{-1}...V_1\rme^{z_1D_1}V_1^{-1}
\right)\nonumber\\
&=&\sum_{r=1}^R\del_r\exp\left(\sum_{s=1}^mz_s\overline{\gam_{s,r}}\right).\label{piecewisetilde}
\end{eqnarray}
The coefficient $\overline{\gam_{t,r}}$ can be determined by
putting
\[
z_s=\left\{\begin{array}{ll} w&\mbox{ if $s=t$}\\ 0&\mbox{
otherwise}
\end{array}\right.
\]
to obtain
\begin{eqnarray*}
\sum_{r=1}^R\del_r\exp\left(w\overline{\gam_{t,r}}\right)
&=&\det\left(S+TV_t\rme^{wD_t}V_t^{-1}\right)\\
&=&\det\left(V_t^{-1}\right)\det\left(SV_t+TV_t\rme^{wD_t}\right).
\end{eqnarray*}
By expanding the final determinant we deduce that each
$\overline{\gam_{t,r}}$ is the sum of some subset of the numbers
$a_{t,j}^{-1}(\alp_t-\alp_{t-1})$ where $1\leq j\leq n$ and
$\{a_{t,j}\}_{j=1}^n$ are the eigenvalues of $A_t$. Since
$F(z)=\tilde F(z,...,z)$ we also have
\[
\gam_r=\sum_{s=1}^m \gam_{s,r}.
\]
We say that we are in the generic case if all possible $\del_r$ are
non-zero, and also no $\gam_r$ lies within an edge of the convex
hull $K$ of $\{\gam_r\}_{r=1}^R$.

We now apply Corollary 7 exactly as in the proof of Theorem
\ref{asymptheorem} but with different choices for $K_s$.
Generically each $K_s$ is the convex hull of all numbers of the
form
\[
v_{s,J}=\sum_{j\in J} a_{s,j}^{-1}
\]
where $J$ ranges over all subsets of $\{1,...,n\}$. Hence
\[
K_s=I_{s,1}+...+I_{s,n}
\]
where $I_{s,j}$ is the interval with end-points $0,a_{s,j}^{-1}$.
It follows by (\ref{convsum}) that
\[
b(K_s)=\sum_{j=1}^n b(I_{s,j})=2\sum_{j=1}^n |a_{s,j}|^{-1}.
\]
This yields the statement of the theorem.

We had expected that the above theorem would be valid for periodic
boundary conditions, but the following example shows that this need
not be the case.

\begin{example}
Let $\alp=0$, $\bet=2$, $n=2$, $A(x)=A_1$ for $0\leq x <1$ and
$A(x)=A_2$ for $1\leq x\leq 2$. Suppose also that
\[
A_i=V_i \left[\begin{array}{cc}
a_{1,i}&0\\0&a_{2,i}\end{array}\right] V_i^{-1}
\]
for $i=1,2$, where all $a_{j,i}$ are non-zero. Finally put
\[
D_i=\left[\begin{array}{cc} u_i&0\\0&v_i\end{array}\right]
\]
for $i=1,2$, where $u_i=a_{1,i}^{-1}$ and $v_i=a_{2,i}^{-1}$.

Periodic boundary conditions correspond to the choice $S=I$ and
$T=-I$ and lead to the formula
\begin{eqnarray*}
F(z)&=&\det\left(I-V_2\rme^{D_2z}V_2^{-1}.V_1e^{D_1z}V_1^{-1}\right)\\
&=&c\det\left( X-\rme^{D_2z}Xe^{D_1z}\right)
\end{eqnarray*}
where $c=\det(V_2)\det(V_1^{-1})\not= 0$ and $X=V_2^{-1}V_1$ is
invertible. Hence
\begin{eqnarray}
F(z)&=&c\det\left[\begin{array}{cc}
X_{1,1}(1-\rme^{(u_2+u_1)z})&X_{1,2}(1-\rme^{(u_2+v_1)z})\\
X_{2,1}(1-\rme^{(v_2+u_1)z})& X_{2,2}(1-\rme^{(v_2+v_1)z})
\end{array}\right]. \label{periodicF}
\end{eqnarray}
If all $X_{i,j}$ are non-zero we deduce that the possible values of
$\gam_r$ are $0,u_1+u_2,u_1+v_2,v_1+u_2,v_1+v_2$ and
$u_1+u_2+v_1+v_2$. Depending on the positions of these points $K$
may have several shapes, and hence $b(K)$ may have several values.
One possibility is
\[
b(K)=2|u_2+u_1|+2|v_2+v_1|
\]
which is quite different from the formula obtained in Theorem
\ref{genericasymp}.

In the very special case $A_2=-A_1$ \emp{every} solution of $Lf=zf$
is periodic and $\Spec(L)=\C$. But in this case $u_2+u_1=v_2+v_1=0$
and $X=I$, so (\ref{periodicF}) yields $F(z)=0$ for all $z$.

\end{example}

\section{Basis Problems}

The fact that one may be able to determine the eigenvalues and
eigenfunctions of an operator does not imply that the
eigenfunctions form a basis in the relevant Banach space,
\cite{ad,dav2}. If this fails then the positions of the
eigenfunctions may be very unstable with respect to small
perturbations of the operator, and the significance of the spectrum
becomes moot. Boulton \cite{bou,bou2} has investigated a closely
related problem for second order elliptic systems in the language
of pseudospectral theory, \cite{bot,dav,dav3,red,rt,tre}, which
amounts to estimating the resolvent norms of the operators
concerned.

In this section we show that such problems do indeed occur in the
context of this paper. For some operators of the type we consider
the eigenfunctions form a basis, but for many others, and we
believe most, they do not. We do not attempt a complete analysis,
but just discuss the simplest example, of a first order system with
two components. The operator $H=L^2$ provides a second order system
with the same eigenvectors as $L$ and hence the same basis
problems.

\begin{example}\label{satwodim}
Let $n=2$ and $0\not= t\in \R$ and let $L$ be defined by
\begin{eqnarray*}
(Lf)_1(x)&:=& tif_1^\pr(x)\\ (Lf)_2(x)&:=&-tif_2^\pr(x)
\end{eqnarray*}
for all $f\in W^{1,2}((0,\pi),\C^2)$, subject to the boundary
conditions $f_1(0)=f_2(0)$ and $f_1(\pi)=f_2(\pi)$. A direct
calculation shows that $\Spec(L)=t\Z$ and that the eigenfunctions
form a complete orthonormal set. Thus $L$ is self-adjoint.
\end{example}

The behaviour of the following slightly more general example is
quite different.

\begin{example}\label{twodim}
Let $n=2$, $s>0$, $t>0$ and $u:=s+it$. Let $L$ be defined by
\begin{eqnarray*}
(Lf)_1(x)&:=& uf_1^\pr(x)\\ (Lf)_2(x)&:=&\overline{u}f_2^\pr(x)
\end{eqnarray*}
for all $f\in W^{1,2}((0,\pi),\C^2)$, subject to the boundary
conditions $f_1(0)=f_2(0)$ and $f_1(\pi)=f_2(\pi)$. A direct
calculation shows that $\lam$ is an eigenvalue if
\[
\rme^{\pi\lam(1/u-1/\overline{u})}=1
\]
or equivalently if
\[
\lam_n=(s^2+t^2) n/t
\]
for some $n\in\Z$. The corresponding eigenfunction is
\begin{eqnarray*}
f_{n,1}(x)&=&\rme^{x\lam_n/u}\\
f_{n,2}(x)&=&\rme^{x\lam_n/\overline{u}}.
\end{eqnarray*}
\end{example}

\begin{lemma} \label{alg}
Let $S$ be a compact subset of $\C$ with zero
Lebesgue measure and suppose that $\C\backslash S$ has two
components, one containing $0$ and the other unbounded. Then the
linear span of the functions $\{z^n\}_{n\in\Z}$ is uniformly dense
in $C(S)$.
\end{lemma}

\Proof By Hartogs-Rosenthal lemma, \cite[p 47]{gam}, $C(S)=R(S)$
where the latter is the uniform closure in $C(S)$ of the space of
rational functions which do not have a pole on $S$. Each such
rational function may be written as a linear combination of
functions $z^m$, $(z-\sig )^{-n}$ where $m,n\geq 0$ and $\sig\notin
S$. If $\sig$ is in the unbounded component of $\C\backslash S$
then $(z-\sig )^{-n}$ may be uniformly approximated by polynomials
in $z$ by Runge's theorem, \cite[p 28]{gam}. If $\sig$ is in the
bounded component then $(z-\sig )^{-n}$ may be uniformly
approximated by polynomials in $z^{-1}$ by using inversion and
Runge's theorem. Putting these facts together completes the proof.

\begin{theorem}
The set of eigenfunctions $\{f_n\}_{n\in\Z}$ is complete in the
sense that its linear span is dense in $L^2((0,\pi),\C^2)$.
\end{theorem}

\Proof Define $w:[0,\pi]\to \C$ by
\[
w(x)=f_{1,1}(x)=\rme^{x(\tau-i)}
\]
where $\tau=s/t$. Let $S=T\cup\overline{T}$ where $T=w([0,\pi])$.
Then $S$ is a closed curve in $\C$ surrounding the origin and
satisfies the conditions of Lemma \ref{alg}.

It is sufficient to prove that if $\phi\in C([0,\pi],\C^2)$
satisfies $\phi(0)=\phi(\pi)=0$ then $\phi$ may be uniformly
approximated by finite linear combinations of $\{f_n\}_{n\in\Z}$.
Given such a $\phi$ we define $\psi\in C(S)$ by
\begin{eqnarray*}
\psi(w(x))&=& \phi_1(x) \mbox{  if $x\in [0,\pi]$}\\
\psi(\overline{w(x)})&=& \phi_2(x) \mbox{  if $x\in [0,\pi]$}.
\end{eqnarray*}
Given $\eps >0$ there exists an approximation
\[
\norm \psi(z)-\sum_{r=-N}^N\del_r z^r \norm_\infty <\eps
\]
by Lemma \ref{alg}. Putting $z=w(x)$ where $x\in [0,\pi]$ we obtain
\[
\norm \phi_1- \sum_{r=-N}^N\del_r f_{r,1} \norm_\infty <\eps.
\]
Putting $z=\overline{w(x)}$ where $x\in [0,\pi]$ we obtain
\[
\norm \phi_2- \sum_{r=-N}^N\del_r f_{r,2} \norm_\infty <\eps.
\]
Combining these we obtain the required estimate
\[
\norm \phi- \sum_{r=-N}^N\del_r f_{r} \norm_\infty <2\eps.
\]

The fact that the eigenvalues of $L$ are real does not imply that
it is similar to a self-adjoint operator.

\begin{theorem}\label{projtheo}
Let $P_n$ be the spectral projection associated with the eigenvalue
$\lam_n$ of $L$. Then $\norm P_n\norm$ diverges at an exponential
rate as $n\to\infty$. The eigenfunctions $f_n$ therefore cannot
constitute a basis of $L^2((0,\pi),\C^2)$.
\end{theorem}

\Proof We have
\[
P_n\phi=\frac{\langle \phi,g_n\rangle f_n}{\langle f_n,g_n\rangle}
\]
where $g_n$ is the appropriate eigenfunction of $L^\ast$, and hence
\[
\norm P_n \norm =\frac{\norm g_n\norm\,\norm f_n \norm}
{|\langle f_n,g_n\rangle |}.
\]
A direct calculation shows that the eigenfunction $g_n$ is given by
\begin{eqnarray*}
g_{n,1}(x)&=&u\rme^{-x\lam_n/\overline{u}}\\
g_{n,2}(x)&=&-\overline{u}\rme^{-x\lam_n/u}.
\end{eqnarray*}
Evaluating the relevant integrals we find that
\begin{eqnarray*}
\langle f_n,g_n\rangle&=& -2\pi it\\
\norm f_n \norm^2 &\sim&\frac{t}{sn}\rme^{2\pi sn/t}\\
\norm  g_n \norm^2 &\sim &\frac{t}{sn}(s^2+t^2)
\end{eqnarray*}
as $n\to +\infty$, with a similar formula as $n\to -\infty$. This
implies the statement about $\norm P_n \norm$. The final statement
is a consequence of the fact that if the eigenfunctions form a
basis then $\norm P_n \norm$ must be a bounded sequence,
\cite{ad,dav2,gk}.

\begin{note}
If $\phi\in L^2$ and $\phi_n=P_n\phi$ then it is possible that
\begin{equation}
\phi\sim\sum_{n\in\Z}\phi_n\label{summ}
\end{equation}
in the sense of some Abel-type summation scheme,
\cite{lidskii,agran,shkal}. In the particular case in which
$\phi(x)=(1,0)$ for all $x\in(0,\pi )$ one finds that $\norm
\phi_n\norm\to\infty$ exponentially fast as $n\to\pm\infty$. It is
not clear that the convergence of (\ref{summ}) using a summation
scheme would have much numerical significance, because of the high
instability of the spectrum under small perturbations of the
operator.
\end{note}

\section{Numerical Range Problems}

The properties which we have discussed so far are similarity
invariants. In other words if $L$ has the property and $T$ is a
bounded invertible operator then $TLT^{-1}$ also has the property.
In this section we discuss properties which depend upon the
particular norm chosen out of a similarity class. These are
important because the norm is often given by physical
considerations, and even if it is not, one might not know whether
some better equivalent norm exists or how to find it.

Let $Hf=-Af^{\pr\pr}$ in $L^2((0,\bet ),\C^n)$ where $A$ is a
bounded invertible $n\times n$ matrix. We impose any linear
boundary condition at $x=\bet$ and a boundary condition at $x=0$ of
the form $f(0)\in U$, $f^\pr (0)\in V$, where $U$, $V$ are linear
subspaces of $\C^n$.

The following theorem, which extends results of Boulton,
\cite{bou,bou2}, establishes that Kato's theory of sectorial forms,
one of the main tools in non-self-adjoint semigroup theory, cannot
be applied to such operators even if the eigenvalues all lie in a
half plane $\{ z\in
\C:\Re(z)\geq k\}$ for some $k$, as one might expect if the
eigenvalues of $A$ all have positive real parts.

\begin{theorem}\label{nr}
If there exist $c\in U$ and $d\in V$ such that $\langle
Ad,c\rangle\not= 0$ then the numerical range of $H$ equals the
entire complex plane.
\end{theorem}

\Proof Let $\phi\in C^\infty (\R)$ satisfy $\phi(x)=1$ if $x\leq
1/3$ and $\phi(x)=0$ if $x\geq 2/3$. Let also $\psi\in\R$ and
$n\in\Z_+$, and define
\[
f_n(x)=\left(
c+\rme^{i\psi}d\left[(x+1/n)^{2/3}-(1/n)^{2/3}\right]\right)\phi(x)
\]
so that $f\in\Dom(H)$. An easy calculation shows that
\[
\lim_{n\to\infty} \norm f_n \norm^2=\int_0^\bet
\left|\left( c+\rme^{i\psi}dx^{2/3}\right)\phi(x)\right|^2\rmd
x
\]
which we denote by $k_1\not= 0$. An integration by parts
establishes that
\[
\langle Hf_n,f_n\rangle =Q(f_n)+B_n
\]
where
\begin{eqnarray*}
Q(f_n)&=&\int_0^\bet \langle Af_n^\pr,f_n^\pr\rangle \rmd x\\ &\to&
\int_0^\bet \left|  \left(
c+\rme^{i\psi}dx^{2/3}\right)\phi^\pr(x)
+\mbox{$\frac{2}{3}$}\,\rme^{i\psi}dx^{-1/3}\phi(x)\right|^2\rmd
x\\
 &= & k_2
\end{eqnarray*}
as $n\to\infty$, and
\[
B_n=\langle
Af_n^\pr(0),f_n(0)\rangle=\mbox{$\frac{2}{3}$}\,\rme^{i\psi}n^{1/3}\langle
Ad,c\rangle.
\]
We deduce that
\[
\frac{\langle Hf_n,f_n\rangle}{\langle
f_n,f_n\rangle}\sim
\frac{k_2+\frac{2}{3}\rme^{i\psi}n^{1/3}\langle Ad,c\rangle}{k_1}
\]
as $n\to\infty$. By varying $\psi\in\R$ and using the fact that the
numerical range is a convex set, \cite[Theorem 6.1]{ops}, we deduce
that it must equal $\C$.

\begin{example} We consider the operator $H$ defined on
$L^2((0,\pi),\C^2)$ by
\begin{eqnarray*}
(Hf)_1&=&-\rme^{2i\alp}f_1^{\pr\pr}\\
(Hf)_2&=&-\rme^{-2i\alp}f_2^{\pr\pr}
\end{eqnarray*}
subject to the four boundary conditions $f_1(0)=0$, $f_2^\pr(0)=0$
and
\begin{eqnarray*}
\cos(\theta)f_1(\pi)+\sin(\theta)f_2(\pi)&=&0\\
-\sin(\theta) f_1^\pr(\pi)+\cos(\theta)f_2^\pr(\pi)&=&0.
\end{eqnarray*}
To avoid redundancy we assume that $-\pi/2<\theta\leq \pi/2$.

In the particular case $\theta=0$ the two components of $H$ are
independent and its eigenvalues are $n^2\rme^{2i\alp}$ where
$n=1,2,...$ and also $m^2\rme^{-2i\alp}$ where $m=0,1,2,...$
Moreover $H$ is normal and its eigenfunctions form a complete
orthonormal set in $L^2$.

If $\theta\not= 0$ then we may apply Theorem \ref{nr} with
$c=(\sin(\theta),-\cos(\theta))$ and
$d=(\cos(\theta),\sin(\theta))$. We obtain
\[
\langle Ad,c\rangle =i\sin(2\theta)\sin(2\alp)
\]
which is generically non-zero. We deduce that the numerical range
of $H$ equals $\C$. This fact cannot, however, be discovered simply
by finding the eigenvalues. If $\theta\not= 0$ then $0$ is not an
eigenvalue of $H$. The eigenfunction associated to the eigenvalue
$\lam=z^2$ must be of the form
\begin{eqnarray*}
f_1(x)&=&c_1\sinh\left(xz\rme^{-i\alp}\right)\\
f_2(x)&=&c_2\cosh\left(xz\rme^{i\alp}\right)
\end{eqnarray*}
if it is to satisfy the first two boundary conditions. The
existence of a non-zero eigenfunction satisfying the other two
boundary conditions then forces $F(z)=0$ where
\begin{eqnarray*}
F(z)&:=&\rme^{i\alp}\cos(\theta)^2
\sinh(\pi z\rme^{-i\alp})
\sinh(\pi z\rme^{i\alp})\\ &&+\,\rme^{-i\alp}\sin(\theta)^2
\cosh(\pi z\rme^{-i\alp})
\cosh(\pi z\rme^{i\alp}).
\end{eqnarray*}
This is of the canonical form
\[
F(z)=\sum_{r=1}^4 \del_r \rme^{z\overline{\gam_r}}
\]
where $\gam_1=2\pi\cos(\alp)$, $\gam_2=2\pi i\sin(\alp)$,
$\gam_3=-2\pi\cos(\alp)$ and $\gam_4=-2\pi i\sin(\alp)$. The points
$\gam_r$ are the vertices of a rhombus, and the directions of the
four outward pointing normals are $\pm i\rme^{\pm i\alp}$. For
$|\alp|<\pi/4$ this implies that all except finitely many of the
eigenvalues $\lam=z^2$ lie in the half plane $\C^-=\{
z:\Re(z)<0\}$. We conjecture that
\emp{all} the eigenvalues of $H$ satisfy $\Re (\lam)<0$.
\end{example}

\section{Semigroup Properties}

In spite of the above, there are cases in which one can prove that
the operator $H$ is the generator of a contraction semigroup on
$L^2((\alp,\bet),\C^n)$ provided this space is given a new
equivalent norm. The following two theorems describe the abstract
situation.

\begin{theorem}\label{si} Let $Z$ be a closed operator on the
Hilbert space $\cH$ and let $A$ be a bounded invertible operator.
If $0\notin \Spec (Z)$ and $AZ$ is accretive then $-AZ$ is the
generator of a contraction semigroup and hence $\Spec(AZ)\subseteq
\C^+$, the set of complex numbers with non-negative real parts.
\end{theorem}

\Proof The operator $AZ$ is closed and invertible, so there exists
$\eps >0$ such that $z\notin\Spec(AZ)$ for all $|z|<\eps$. The
theorem now follows by applying \cite[Theorem 2.25]{ops}.

\begin{theorem}\label{co} Let $Z$ be a closed operator on the
Hilbert space $\cH$ and let $A$ be a bounded invertible operator.
If $Z+\lam I$ is invertible for some $\lam >0$ and $AZ,\,A$ are
both accretive, then $-AZ$ is the generator of a contraction
semigroup and hence $\Spec(AZ)\subseteq
\C^+$.
\end{theorem}

\Proof By the hypotheses the operator $A(Z+\lam I)$ is invertible
and accretive. Hence $-A(Z+\lam I)$ generates a contraction
semigroup by Theorem \ref{si}. Since $A$ is a bounded perturbation,
$-AZ$ generates a strongly continuous semigroup. This must be a
contraction semigroup by \cite[Theorem 2.27]{ops}.

We now return to a more concrete context. Let $ Z=-\rmd^2/\rmd x^2$
act in $L^2((\alp,\bet),\C^n)$ subject to the boundary conditions
$f(\alp )\in U$, $f(\bet )\in U$, $f^\pr (\alp )\in V$, $f^\pr
(\bet )\in V$ where $\dim(U)+\dim(V)=n$.

\begin{theorem} If $U\cap V\not= \{0\}$ then $\Spec (Z)=\C$. If
$U\cap V= \{0\}$ then $Z$ is similar to a self-adjoint operator
with non-negative discrete spectrum.
\end{theorem}

\Proof Let $0\not= c\in (U\cup V)^\perp$. Then for any $z\in\C$ the
function
\[
g(x)=c\rme^{i\overline{z}x}
\]
is orthogonal to the range of $Z-z^2 I$, so $z^2\in\Spec(Z)$.

If, on the other hand, $U\cap V= \{0\}$ then we may write $Z$ as
the (non-orthogonal) direct sum of its restrictions to
$L^2((\alp,\bet),U)$ and $L^2((\alp,\bet),V)$. In the first we are
imposing Neumann boundary conditions, and in the second Dirichlet
boundary conditions, so the spectrum is as stated. Another way to
express the same idea is  to give $\C^n$ an inner product which
makes $U$ and $V$ orthogonal, so that $Z$ becomes self-adjoint with
respect to this inner product. The change from one inner product to
an equivalent one amounts to the same as a similarity
transformation.

\begin{note}
In the above theorem $0\in\Spec(Z)$ unless $U=\{0\}$ and
$V=\C^n$, the full Dirichlet case. This is why we need both of
Theorems \ref{si} and \ref{co}.
\end{note}

\begin{note}
If the angle between $U$ and $V$ is very small then the condition
number of the similarity transformation is large and it will be
difficult to distinguish between the two alternative conclusions of
the above theorem. In other words $Z$ will have bad pseudospectral
properties.
\end{note}

The following provides a partial converse to Theorem \ref{nr}. We
have not been able to find a similar result if the boundary
conditions are different at the two ends of the interval.

\begin{theorem}
Let $Hf=-Af^{\pr\pr}$ in $L^2((\alp,\bet),\C^n)$ where $A$ is a
bounded invertible $n\times n$ matrix. We impose boundary
conditions $f(\alp )\in U$, $f(\bet )\in U$, $f^\pr (\alp )\in V$,
$f^\pr(\bet )\in V$ where $U\cap V=\{0\}$ and $U+V=\C^n$. Then if
$A$ is accretive and $AV\perp U$, the operator $-H$ generates a
contraction semigroup in $L^2((\alp,\bet),\C^n)$.
\end{theorem}

\Proof We can apply Theorem \ref{co} with $Zf=-f^{\pr\pr}$ subject
to the above boundary conditions provided we show that $H$ is
accretive. If $f\in\Dom(H)$ then integration by parts yields
\[
\langle Hf,f\rangle =
\int_\alp^\bet \langle Af^\pr,f^\pr\rangle \rmd x.
\]
This has non-negative real part since $A$ is accretive.

The condition $AV\perp U$ depends upon the particular inner product
used in $\C^n$ and is not satisfied generically. On the other hand
the condition $AV\cap U=\{0\}$ does hold generically. The second
condition implies the first for a suitable choice of the inner
product on $\C^n$, but one then needs to check whether $A$ is
accretive for this new inner product. Assuming this does happen, if
the condition number relating the two inner products is $\kap$ then
one only gets
\[
\norm \rme^{-Ht}\norm\leq\kap
\]
with respect to the standard inner product for all $t\geq 0$. This
may be of value if $\kap$ is fairly small but it is of little
computational use if $\kap$ is sufficiently large.

\begin{example}\label{periodic} One can prove that $L$ generates
a one-parameter group in the follow situation. Let $Lf=Af^\pr$
subject to the quasi-periodic boundary conditions $f(\bet )=Sf(\alp
)$, where $A$ is a diagonal matrix with real, non-zero eigenvalues
and $S$ is invertible. It appears that $L$ generates a periodic
flow on $[\alp,\bet]$, but with different flow rates in the
different components and mixing at the end points. One can also
regard this as a flow in a generally non-trivial vector bundle over
the circle. We restrict attention to the simple case $A=I$, but
note that Boulton's ideas allow one to treat the case in which the
eigenvalues of $A$ are real and rationally related,
\cite{bou,bou2}.

Under the above assumptions we have
\[
F(z)=\det\left\{ S-\rme^{z(\bet-\alp )}I\right\}.
\]
Putting $\sig=\rme^{z(\bet-\alp )}$ and finding the eigenvalues of
$S$ we obtain $n$ generally distinct solutions $\sig=\sig_r$. We
conclude that the eigenvalues of $L$ lie on one of $n$ distinct
lines parallel to the $y$-axis.
\end{example}

Our next theorem implies that it is only possible to prove
that $\rme^{Lt}$ is a strongly continuous one-parameter group using
numerical range ideas if $L$ is already skew-adjoint. In spite of the
failure of this method, the subsequent theorem shows that $L$ does
generate a one-parameter group whenever the eigenvalues of $A$ are
all real.

\begin{theorem}
The numerical range of $L$ is only contained in a strip of the form
\[
\{ z\in\C:c_1\leq \Re (z)\leq c_2\}
\]
if $L$ is already skew adjoint, in other words if $S$ is unitary.
\end{theorem}

\Proof  For general functions in $\Dom(L)$ we have
\begin{eqnarray*}
\Re \langle Lf,f\rangle &=& \norm f(\bet )\norm^2-\norm f(\alp )\norm^2 \\
&=& \langle (S^\ast\, S-I)f(\alp ),f(\alp )\rangle.
\end{eqnarray*}
If $S$ is not unitary then this may be arbitrarily large (positive
or negative or both) for functions in $\Dom(L)$ of unit norm; the
proof is similar to that of Theorem \ref{nr}.  This prevents the
numerical range from lying in a strip of the stated kind.

\begin{theorem}
For all operators $L$ defined as in Example \ref{periodic}, $L$ is
the generator of a strongly continuous one-parameter group acting
on $L^2(\alp,\bet),\C^n)$.
\end{theorem}

\Proof We define the bounded invertible operator $W:L^2\to L^2$ by
\[
W f(x)=S(x)f(x)
\]
where $S(x)$ is any smooth function with values in the invertible
$n\times n$ matrices satisfying $S(\alp )=I$ and $S(\bet )=S$. Now
let $\cD$ be the space of functions in $W^{1,2}$ such that $f(\alp
)=f(\bet )$ and define $M:\cD\to L^2$ by
\[
Mf=W^{-1}L W f.
\]
A direct calculation shows that
\[
Mf(x)=f^\pr (x)+Y(x)f(x)
\]
where $Y$ is a bounded smooth matrix-valued function on
$[\alp,\bet]$. Writing this is the more abstract form $M=M_0+Y$ we
see that $M$ is a bounded perturbation of $M_0$, which generates a
one-parameter group of isometries on $L^2$. Hence $M$, and then
$L$, generate one-parameter groups on $L^2$.

\section{Inverse Spectral Theory}

In \cite{bou,bou2}, Boulton showed that a second order elliptic
system $H$ may have real spectrum without being similar to a
self-adjoint operator. In this section we investigate the extent to
which one can reconstruct a first order elliptic system from its
spectrum. We start with a partial positive result and then give an
example which shows that in general one cannot go any further.

We continue in our earlier framework including the generic
assumption of Section~3. The first condition of the following
theorem is probably unnecessary. We conjecture that in general the
multiplicity of any eigenvalue of $L$ equals the order of the
corresponding zero of $F(z)$. Note that we have already proved that
all zeros of $F(z)$ of large enough modulus are simple, so the
condition only concerns a finite number of smaller eigenvalues.

\begin{theorem}\label{polygon}
Assume that the zeros of $F(z)$ are all simple. Then the spectrum
of $L$ determines the function $F(z)$ up to an exponential factor,
and hence the set of constants $\gam_r$, $1\leq r\leq R$, up to a
common additive factor.
\end{theorem}

\Proof Although we do not use this in the proof, we note that
by Theorems \ref{gr} and \ref{f} the
asymptotics of the spectrum of $L$ determine the length and
direction of each edge of $K$. Since $K$ is a convex polygon it can
only be reconstructed in one way up to translations in $\C$.

The statement of the theorem is a consequence of the fact that
$F(z)$ is an entire function of order $1$, and can therefore be
written in the form
\begin{equation}
 F(z)=z^m\rme^{hz}\prod_{n=1}^\infty \left\{
\left(1-\frac{z}{z_n}\right)\rme^{z/z_n}\right\}\label{product}
\end{equation}
for some $h\in\C$, where $z_n$ are the zeros of $F(z)$. See
\cite[p. 199]{hille}. The fact that $F(z)$ determines all of the
$\gam_r$, $1\leq r\leq R$, up to the additive factor $\overline{h}$
is elementary.

\begin{note}
One may break the product in (\ref{product}) into sporadic terms
together with terms associated with $z_n$ which are close to one of
the straight lines already described. The products associated with
each line are asymptotically similar to certain formulae involving
gamma functions, namely
\[
\frac{\rme^{(\gam+\del) z}\Gam(a+1)}{\Gam(a-z+1)}=
\prod_{n=1}^\infty \left\{
\left(1-\frac{z}{a+n}\right)\rme^{z/(a+n)}\right\}
\]
where $\gam$ is Euler's constant and
\[
\del=\sum_{n=1}^\infty \frac{a}{(a+n)n}.
\]
This follows immediately from Whittaker and Watson, \cite[p
236]{ww}.
\end{note}

In certain cases one can go further than Theorem \ref{polygon}. If
$\dim(W)=1$ then we showed in Section 2 that $\gam_r=(\bet-\alp
)/\overline{a_r}$ for $1\leq r\leq R$ where $a_r$ are the
eigenvalues of $A$. If in addition we know the centre of $K$ and
the value of $\bet-\alp$, then we can evaluate all of the
eigenvalues of $A$. The following example shows that in general one
cannot go any further.

\begin{example}\label{counter}
Given $k>0$ let $u=s+it$ lie on the circle $\cC=\{u: s^2+t^2=k t\}$
and consider again Example \ref{twodim}. We exclude the case
$s=t=0$ for obvious reasons; the case $s=0,\, t=\sqrt{k}$ is
covered by Example \ref{satwodim}. The two vertices of $K$ are
$\gam_\pm=\pi(v\pm ik)$ where $v=s/(s^2+t^2)=s/kt$. We thus see
that $K$ is the same, up to a translation, for all of these
operators as $u$ varies on $\cC$. By Examples \ref{satwodim} and
\ref{twodim} the spectrum of $L$ is $k\Z$ for all choices of $u$ on
$\cC$.

In order to see that $L$ cannot be reconstructed from its spectrum
up to a similarity transformation, we need to prove the following
theorem.
\end{example}

\begin{theorem}
No two of the operators of Example \ref{counter} are similar as $u$
varies in $\cC$.
\end{theorem}

\Proof We use the notation of Theorem \ref{projtheo}. If two such
operators $L$ and $\tilde{L}$ determined by $u,\tilde{u}\in\cC$
were similar then it would follow that $\norm
\tilde{P_n}\norm/\norm P_n\norm$ would be bounded above and below
by positive constants uniformly with respect to $n$.

Now the proof of Theorem \ref{projtheo} implies that
\[
\norm P_n\norm\sim \frac{\rme^{\pi k n v}|u|}{2\pi |t| k n v }
\]
as $n\to\infty$. We deduce that if $L$ and $\tilde{L}$ are similar
then $v=\tilde{v}$. This implies that $u=\tilde{u}$.

\vskip 1in
%%%%%%%%%%%%%%%%%%%%%%%%%
{\bf Acknowledgments } I would like to thank J~M~Anderson,
L~S~Boulton, K~M~Ball, D~A~Edwards, G~Rozenblum, S~G~Scott and
R~F~Streater for valuable discussions.
\vskip 0.5in
%\newpage
%%%%%%%%%%%%%%%%%%%%%%%

\par
%%%%%%%%%%%%%%%%%%%%%%%%%%%%%%%%%%%%%%%%
\vskip 0.5in
Department of Mathematics \newline King's College \newline Strand
\newline London WC2R 2LS \newline England\\
\vskip 0.1in e-mail: E.Brian.Davies@kcl.ac.uk
\vfil
\end{document}